\documentclass[11pt]{amsart}

\usepackage[hidelinks]{hyperref}

\usepackage{amssymb}
\usepackage{amsmath}
\usepackage{mathrsfs}
\usepackage{url}
\usepackage{verbatim}
\usepackage{multicol}
\input xy
\xyoption{all}

\usepackage{graphicx}  
\usepackage{tikz}
\usetikzlibrary{decorations.pathreplacing}
\usepackage{tkz-graph}
\usetikzlibrary{arrows}
\pgfarrowsdeclare{mytip}{mytip}{%
  \setlength{\arrowsize}{1pt}
  \addtolength{\arrowsize}{.5\pgflinewidth}
  \pgfarrowsrightextend{0}
  \pgfarrowsleftextend{-5\arrowsize}
}{%
  \setlength{\arrowsize}{1pt}
  \addtolength{\arrowsize}{.5\pgflinewidth}
  \pgfpathmoveto{\pgfpoint{-5\arrowsize}{1.5\arrowsize}}
  \pgfpathlineto{\pgfpointorigin}
  \pgfpathlineto{\pgfpoint{-5\arrowsize}{-1.5\arrowsize}}
  \pgfusepathqstroke
}

\newtheorem{thm}{Theorem}
\newtheorem{lemma}[thm]{Lemma}

\newtheorem{cor}[thm]{Corollary}
\newtheorem{prop}[thm]{Proposition}

\newtheorem{Definition}[thm]{Definition}
\newenvironment{definition}
  {\begin{Definition}\rm}{\end{Definition}}

\newtheorem{Example}[thm]{Example}
\newenvironment{example}
  {\begin{Example}\rm}{\end{Example}}
  
\newtheorem{Remark}[thm]{Remark}
\newenvironment{remark}
  {\begin{Remark}\rm}{\end{Remark}}

\newcommand{\R}{\mathbb{R}}

\usepackage{multirow}
\usepackage{mathdots}
\usepackage{xhfill}

\newcommand*{\bigchi}{\mbox{\Large$\chi$}}
\newcommand{\sbullet}{\scalebox{0.5}{$\bullet$}}

\title{Bigraphical arrangements}
\author{Sam Hopkins}
\email{samuelfhopkins@gmail.com}
\author{David Perkinson}
\email{davidp@reed.edu}
\address{Reed College, Portland OR, 97202}

\begin{document}

\begin{abstract}
We define the bigraphical arrangement of a graph and show that the Pak-Stanley labels of its regions are the parking functions of a closely related graph, thus proving conjectures of Duval, Klivans, and Martin~\cite{duval} and of Hopkins and Perkinson~\cite{hopkins}. A consequence is a new proof of a bijection between labeled graphs and regions of the Shi arrangement first given by Stanley in \cite{pak-stanley}. We also give bounds on the number of regions of a bigraphical arrangement.
\end{abstract}

\maketitle


Throughout this paper $G$ will be a simple graph (no multiedges or loops, but not necessarily connected) with vertex set $V = \{v_1,\ldots,v_n\}$ and edge set $E$. A \emph{(real) hyperplane arrangement} is a finite collection of affine hyperplanes in Euclidean space. Our object of study is the \emph{bigraphical arrangement} of $G$, so called because it associates two hyperplanes to each edge of the graph. 

\begin{definition}
For each edge $\{v_i,v_j\} \in E$, choose real numbers $a_{ij}$ and $a_{ji}$ such that there exists $x \in \mathbb{R}^n$ with $x_i - x_j < a_{ij}$ and $x_j - x_i < a_{ji}$ for all $\{v_i,v_j\} \in E$. We call these numbers \emph{parameters} and we call $A := \{a_{ij}\}$ a \emph{parameter list}. The \emph{bigraphical arrangement} $\Sigma_G(A)$ is the set of $2|E|$ hyperplanes,
\[ \Sigma_G(A) := \{x_i - x_j = a_{ij}\colon \{v_i,v_j\} \in E\}. \]
The \emph{regions} of $\Sigma_G(A)$ are the connected components of $\mathbb{R}^{n} \setminus \Sigma_G(A)$. The \emph{central region} is the region defined by $x_i - x_j < a_{ij}$ for all $\{v_i,v_j\} \in E$. The above condition on $A$ guarantees that the central region is nonempty.
\end{definition}

Several authors have connected hyperplane arrangements to graphs in various ways. The \emph{graphical arrangement}~\cite[p. 414]{stanley} of $G$,
\[ \{x_i - x_j = 0\colon \{v_i,v_j\} \in E \textrm{ with $i<j$}\},\]
associates a single hyperplane to each edge of $G$. One interesting property of the graphical arrangement is that its characteristic polynomial is the chromatic polynomial of $G$. Thus, the graphical arrangement encodes information about colorings of $G$. As we shall see, the bigraphical arrangement encodes information about the parking functions of $G_{\sbullet}$, the graph obtained from $G$ by adding a special sink vertex $v_0$ and an edge between $v_0$ and each vertex $v \in V$. For background on parking functions, see \S2.

J.-Y. Shi~\cite{shi}, in his investigation of the Kazhdan-Lusztig cells of affine Weyl groups of type $A_{n-1}$, introduced the \emph{Shi arrangement},
\[ \{x_i - x_j = 0,1\colon 1 \leq i < j \leq n\}. \]
He proved that the number of regions of this arrangement is $(n+1)^{n-1}$, Cayley's formula for the number of spanning trees of the complete graph~$K_{n+1}$. Stanley~\cite{pak-stanley}, in collaboration with Pak, was the first to give a bijective proof of this result by labeling the regions of the Shi arrangement with parking functions. (There are several well-known bijections between parking functions and spanning trees.) Stanley and Pak's procedure labels the central region of the Shi arrangement with the parking function $00\ldots0$. It then inductively labels the other regions by moving outwards and increasing the~$i$th coordinate of a region's label whenever a hyperplane is crossed that corresponds to an increase in $x_i$. We call the resulting labels the \emph{Pak-Stanley labels} of the regions of an arrangement. Duval, Klivans, and Martin~\cite{duval} defined the \emph{$G$-Shi arrangement},
\[ \{x_i - x_j = 0,1\colon  \{v_i,v_j\} \in E \textrm{ with $i<j$} \},\]
and conjectured that the Pak-Stanley labels of the $G$-Shi arrangement are the parking functions of $G_{\sbullet}$.\footnote{The $G$-Shi arrangement is not to be confused with what Armstrong and Rhoades~\cite{armstrong} call the \emph{deleted Shi arrangement} and denote $\mathrm{Shi}(G)$:
\[\mathrm{Shi}(G) := \{ x_i - x_j = 0\colon 1 \leq i < j \leq n\} \cup \{ x_i - x_j = 1\colon  \{v_i,v_j\} \in E \textrm{ with $i<j$} \}.\] While the $G$-Shi arrangement has $2|E|$ hyperplanes and is a bigraphical arrangement, the $\mathrm{Shi}(G)$ arrangement has $n + |E|$ hyperplanes and is therefore not a bigraphical arrangement.}  We prove this conjecture as a consequence of Corollary~\ref{cor:main}. The $G$-Shi arrangement is in fact a special kind of bigraphical arrangement. 

In~\cite{hopkins}, Hopkins and Perkinson studied the \emph{$G$-semiorder arrangement},
\[ \{x_i - x_j = 1\colon \{v_i,v_j\} \in E\}, \]
another special kind of bigraphical arrangement. They showed that the Pak-Stanley labels of the $G$-semiorder are the $G_{\sbullet}$-parking functions sought by Duval, Klivans, and Martin. It was also conjectured in~\cite{hopkins} that if one were to slide the hyperplanes of the $G$-semiorder arrangement along their normals, although some regions are destroyed and others are created, so long as the central region is preserved the set of parking function labels remains the same. In this way, one could deform the $G$-semiorder arrangement into the $G$-Shi arrangement and show that the $G$-Shi arrangement has the expected set of labels. Figure~\ref{fig:sliding} depicts this sliding procedure when $G = K_3$. 

Our Corollary~\ref{cor:main} establishes the sliding conjecture: the Pak-Stanley
labels of any bigraphical arrangement, $\Sigma_G(A)$, are the parking
functions of $G_{\sbullet}$. (Remark~\ref{remark:parkfns} indicates how, in
addition, all the parking functions of $G$ with respect to each of its vertices
are encoded in the regions of $\Sigma_G(A)$.)  In proving
Corollary~\ref{cor:main}, we generalize a result of Benson, Chakbarty, and
Tetali~\cite{benson}, who show that acyclic total orientations of $G$ correspond
to maximal parking functions of $G_{\sbullet}$. We show that certain families of
partial orientations of $G$ defined in \S\ref{sec:pos} correspond to all of the
parking functions of $G_{\sbullet}$. If $G$ is the complete graph $K_n$ and the
parameter list $A$ corresponds to the Shi arrangement, Corollary~\ref{cor:main}
provides an alternate proof of the bijection of Pak and Stanley between regions
of the Shi arrangement and parking functions.  

\begin{figure}
\begin{tikzpicture}[scale=1]
\draw [style=thick] (-3,4) -- (0.5,-4);
\node at (-3.5,4.2) {$x_1-x_3=1$};
\draw [style=thick] (-0.5,4) -- (3,-4);
\node at (-1.2,4.2) {$x_3-x_1=1$};

\draw [->,style=thick] (0.9,0.8) -- (0.3,0.45);
\draw [->,style=thick] (-0.24,3.4) -- (-0.84,3.05);
\draw [->,style=thick] (2.4,-2.6) -- (1.8,-2.95);

\draw [style=thick] (3,4) -- (-0.5,-4);
\node at (3.4,4.2) {$x_2-x_1=1$};
\draw [style=thick] (0.5,4) -- (-3,-4);
\node at (1.0,4.2) {$x_1-x_2=1$};

\draw [->,style=thick] (0.9,-0.8) -- (0.3,-0.45);
\draw [->,style=thick] (-0.24,-3.4) -- (-0.84,-3.05);
\draw [->,style=thick] (2.4,2.6) -- (1.8,2.95);

\draw [style=thick] (-4,1.5) -- (4,1.5);
\node at (3.4,1.8) {$x_3-x_2=1$};
\draw [style=thick] (-4,-1.5) -- (4,-1.5);
\node at (3.4,-1.8) {$x_2-x_3=1$};

\draw [->,style=thick] (0,1.5) -- (0,0.75);
\draw [->,style=thick] (3,1.5) -- (3,0.75);
\draw [->,style=thick] (-3,1.5) -- (-3,0.75);

\node at (0,0) {$000$};

\node at (0,1.9) {$001$};
\node at (0,-1.9) {$010$};
\node at (1.3,1.1) {$001$};
\node at (-1.3,1.1) {$100$};
\node at (-1.3,-1.1) {$100$};
\node at (1.3,-1.1) {$010$};

\node at (1.2,2.5) {$002$};
\node at (-1.2,2.5) {$101$};
\node at (-1.2,-2.5) {$110$};
\node at (1.2,-2.5) {$020$};
\node at (2.5,0) {$011$};
\node at (-2.5,0) {$200$};

\node at (3.4,-2.6) {$021$};
\node at (3.4,2.6) {$012$};
\node at (-3.4,2.6) {$201$};
\node at (-3.4,-2.6) {$210$};
\node at (0,3.8) {$102$};
\node at (0,-3.8) {$120$};

\end{tikzpicture}

\xhrulefill{gray}{4pt}
\vspace{0.3cm}

\begin{tikzpicture}[scale=1]
\draw [style=thick]  (-3,4) -- (0.5,-4);
\node at (-3.5,4.2) {$x_1-x_3=1$};
\draw [style=thick]  (-0.5,4) -- (3,-4);
\node at (-1.2,4.2) {$x_3-x_1=0$};

\draw [style=thick]  (3,4) -- (-0.5,-4);
\node at (3.5,4.2) {$x_2-x_1=0$};
\draw [style=thick]  (0.5,4) -- (-3,-4);
\node at (1.0,4.2) {$x_1-x_2=1$};

\draw [style=thick]  (-4,0) -- (4,0);
\node at (3.5,0.3) {$x_3-x_2=0$};
\draw [style=thick]  (-4,-2.84) -- (4,-2.84);
\node at (3.5,-2.54) {$x_2-x_3=1$};

\node at (0,-1) {$000$};

\node at (0,1) {$001$};
\node at (1.25,-2) {$010$};
\node at (-1.25,-2) {$100$};

\node at (3,-1) {$011$};
\node at (-3,-1) {$200$};
\node at (1.5,-3.5) {$020$};
\node at (-1.5,-3.5) {$110$};
\node at (-1.2,2.5) {$101$};
\node at (1.2,2.5) {$002$};

\node at (-3.3,-3.25) {$210$};
\node at (3.3,-3.25) {$021$};
\node at (3.3,1.8) {$012$};
\node at (-3.3,1.8) {$201$};
\node at (0,3.8) {$102$};
\node at (0,-3.8) {$120$};

\end{tikzpicture}

\caption{Hyperplanes in the $K_3$-semiorder arrangement (above) are slid in the directions indicated by the arrows to yield the $K_3$-Shi arrangement (below). The set of Pak-Stanley labels, strings~$c_1c_2c_3$ inside each region, remains the same.}\label{fig:sliding}
\end{figure}

\begin{example}
The following are examples of bigraphical arrangements:
\begin{enumerate}
\item Setting $a_{ij} = 1$ for all $i,j$ gives the \emph{$G$-semiorder} arrangement, studied in~\cite{hopkins}. We will use $\mathrm{SEMI}$ to denote the parameter list of the $G$-semiorder arrangement and thus denote the $G$-semiorder arrangement by $\Sigma_G(\mathrm{SEMI})$.
\item Setting $a_{ij} = 1$ if $i < j$ and $0$ otherwise gives the \emph{$G$-Shi} arrangement, the subject of a conjecture in \cite{duval} that we establish as a consequence of Corollary~\ref{cor:main}. We will use $\mathrm{SHI}$ to denote the parameter list of the $G$-Shi arrangement and thus denote the $G$-Shi arrangement by~$\Sigma_G(\mathrm{SHI})$.
\item Let $\eta = (\ell_1,\ldots,\ell_n) \in \mathbb{Z}_{>0}^n$. Setting $a_{ij} = \ell_i$ for all $i,j$ gives what we call the \emph{$(G,\eta)$-interval order arrangement}.
\end{enumerate}
Taking $G$ to be the complete graph $K_n$ recovers the normal semiorder, Shi, and interval order arrangements. See~\cite{stanley} for definitions of these arrangements, as well as for basic concepts from the theory of hyperplane arrangements, in particular, that of the \emph{characteristic polynomial}. $\square$
\end{example}

From now on we assume we have fixed some parameter list $A$. Note that~$\Sigma_G(A)$ having a nonempty central region is essentially equivalent to~$A$ having only positive entries. If the $a_{ij}$ are all positive, then the origin satisfies~$x_i - x_j < a_{ij}$ for all~$\{v_i,v_j\} \in E$. On the other hand, suppose~$\Sigma_G(A)$ has a nonempty central region and that $p$ is a point in this region. Then the translation $x \mapsto x - p$ maps $\Sigma_G(A)$ to a bigraphical arrangement whose parameter list has positive entries.

In \S\ref{sec:pos}, we develop a correspondence between regions of $\Sigma_G(A)$ and partial orientations of $G$. In \S\ref{sec:pfs}, we prove our main result, Corollary~\ref{cor:main}, which says that the Pak-Stanley labeling of any $\Sigma_G(A)$ yields the set of parking functions of $G_{\sbullet}$. In \S\ref{sec:numregs}, we bound the number of regions of $\Sigma_G(A)$ for arbitrary~$A$ and we find its characteristic polynomial when $A$ is generic. The characteristic polynomial of a generic $\Sigma_G(A)$ turns out to be related to the Tutte polynomial of $G$.

\medskip

\noindent {\bf Acknowledgements}: We thank Art Duval, Caroline Klivans, and Jeremy Martin for introducing us to the $G$-Shi conjecture and for helpful comments. We thank students who participated in the geometry class taught by the second author at Reed College during the fall semester of 2012: Max Carpenter, Nadir Hajouji, Austin Humphrey, Mikhail Lepilov, Alex Perusse, Yotam Reshef, Marcus Robinson, Darko Trifunovski,  Chris Vittal, and Austin Young. We also thank Collin Perkinson for help with proofreading.

\section{\texorpdfstring{Partial orientations and the regions of $\Sigma_G(A)$}{Partial orientations and the regions of SG(A)}} \label{sec:pos}

\begin{definition}
A {\em partial orientation} of $G$ is a choice of directions for a subset of the edges of $G$.  Formally, a \emph{step} is an ordered pair $(u,v) \in V\times V$ such that~\mbox{$\{u,v\} \in E$}, and a partial orientation $\mathcal{O}$ is a set of steps with the property that if $(u,v)\in\mathcal{O}$, then $(v,u)\notin\mathcal{O}$. We say $\mathcal{O}$ is \emph{acyclic} if it does not contain a cycle of steps.  \end{definition}

\begin{definition}
Let~$\mathcal{O}$ be a partial orientation.  If $e=\{u,v\}\in E$ and $(u,v)\in\mathcal{O}$, then despite the ambiguity, we write $e\in\mathcal{O}$ and say $e$ is \emph{oriented}.  In that case, we think of~$e$ as an arrow from $u$ to $v$ and write~$e^{-}=u$ and~$e^{+}=v$.  If neither $(u,v)$ nor $(v,u)$ is in $\mathcal{O}$, we write $e\notin\mathcal{O}$ and say that $e$ is an unoriented or {\em blank} edge.  The {\em indegree} of $u\in V$ relative to~$\mathcal{O}$, denoted~$\mathrm{indeg}_{\mathcal{O}}(u)$, is the number of edges~$e\in\mathcal{O}$ such that $e^{+}=u$.  Similarly, the {\em outdegree} of the vertex $u\in V$ relative to $\mathcal{O}$ is the number of edges $e\in\mathcal{O}$ such that $e^{-}=u$.  The \emph{degree} of $u$ is the number of $e\in E$ containing $u$.
\end{definition}

\noindent {\bf Notation.} Partial orientations naturally serve as labels for the regions of bigraphical arrangements. Suppose $R$ is a region of $\Sigma_G(A)$. Define $\mathcal{O}_R$ to be the partial orientation obtained by (i) $(v_i,v_j)\in \mathcal{O}_R$ if $\{v_i,v_j\}\in E$ and~$x_j > x_i + a_{ji}$ in $R$, and~(ii) all other edges are blank. We now classify exactly which partial orientations are labels of regions of $\Sigma_G(A)$.

\begin{definition}
Let $\mathcal{O}$ be a partial orientation. A step $(u,v)$ is \emph{compatible} with $\mathcal{O}$ if $(v,u) \notin \mathcal{O}$. In other words, a step $e$ is compatible with~$\mathcal{O}$ if~$e \in \mathcal{O}$ or $e$ is a blank edge of $\mathcal{O}$.  A \emph{potential cycle for $\mathcal{O}$} is a set $C = \{(u_1,u_2),(u_2,u_3),\ldots,(u_{k},u_{1})\}$ of steps compatible with $\mathcal{O}$.
\end{definition}

\begin{figure}[htb]
\begin{tikzpicture}

\SetVertexMath
\GraphInit[vstyle=Art]
\SetUpVertex[MinSize=3pt]
\SetVertexLabel
\tikzset{VertexStyle/.style = {shape = circle,shading = ball,ball color = black,inner sep = 2pt}}
\SetUpEdge[color=black]

\Vertex[LabelOut,Lpos=90,Ldist=.05cm,x=1,y=4]{v_1}
\Vertex[LabelOut,Lpos=180,Ldist=.05cm,x=0,y=2.5]{v_2}
\Vertex[LabelOut,Lpos=230,Ldist=.05cm,x=0,y=1]{v_3}
\Vertex[LabelOut,Lpos=300,Ldist=.05cm,x=2,y=1]{v_4}
\Vertex[LabelOut,Lpos=0,Ldist=.05cm,x=2,y=2.5]{v_5}
\Edges[style={->,>=mytip,thick}](v_1,v_2)
\Edges(v_1,v_5)
\Edges(v_4,v_3)
\Edges(v_4,v_5)
\Edges[style={->,>=mytip,thick}](v_5,v_2)
\Edges[style={->,>=mytip,thick}](v_2,v_3)

\end{tikzpicture}
\caption{A partial orientation $\mathcal{O}$. Consider the potential cycle $C := \{(v_5,v_2),(v_2,v_3),(v_3,v_4),(v_4,v_5)\}$ for $\mathcal{O}$. We have $\nu_{\mathrm{SEMI}}(C,\mathcal{O}) = 0$, so $\mathcal{O}$ is not $\mathrm{SEMI}$-admissible. However, $\mathcal{O}$ is $\mathrm{SHI}$-admissible.} \label{fig:poex}
\end{figure}
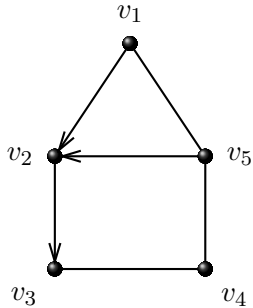

\begin{definition}
Let $\mathcal{O}$ be a partial orientation. The \emph{score} of a step $e = (v_i,v_j)$ compatible with $\mathcal{O}$ is
\[\nu_{A}(e,\mathcal{O}) = \begin{cases}
  \hfill a_{ij} &\text{if $\{v_i,v_j\}\notin \mathcal{O}$},\\
 \hfill -a_{ji}&\text{if $(v_i,v_j) \in \mathcal{O}$}.
\end{cases} \]
The \emph{score} of a potential cycle $C$ for $\mathcal{O}$ is
\[\nu_{A}(C,\mathcal{O}) = \sum_{e \in C} \nu(e,\mathcal{O}).\]
When the parameter list is clear from context, we omit the subscript and we write~$\nu(C,\mathcal{O}) := \nu_{A}(C,\mathcal{O})$.
\end{definition}

\begin{definition}
Let $\mathcal{O}$ be a partial orientation. A potential cycle for $\mathcal{O}$ is \emph{bad} if it has a nonpositive score. We say $\mathcal{O}$ is \emph{$A$-admissible} if it has no bad potential cycles. 
\end{definition}

Figure~\ref{fig:poex} gives an example of a partial orientation that is $\mathrm{SHI}$-admissible but not $\mathrm{SEMI}$-admissible.

\begin{thm} \label{thm:farkas}
The regions of $\Sigma_G(A)$ are in bijection with the $A$-admissible partial orientations of $G$. The bijection is given by $R \mapsto \mathcal{O}_R$.
\end{thm}

\emph{Proof}: Let $\mathcal{O}$ be an $A$-admissible partial orientation, and define $R$ to be the region of $\Sigma_G(A)$ determined by the following inequalities: for each edge~$e = \{v_i,v_j\}$ of $G$:
\begin{itemize}
\item $x_i - x_j < a_{ij}$ and $x_j - x_i < a_{ji}$ if $e \notin \mathcal{O}$;
\item $x_i - x_j < -a_{ji}$ if $(v_i,v_j) \in \mathcal{O}$.
\end{itemize}
We must show that $R$ is nonempty. Encode the system of inequalities defining $R$ in a $k\times n$ real matrix $M$ and a column vector $b$ in $\R^{k}$ such that $x \in R$ if and only if $Mx < b$. We take the rows of $M$ to correspond with steps of~$G$ compatible with~$\mathcal{O}$: a row of $M$ with $1$ in the $i$th entry and $-1$ in the~$j$th entry corresponds to a step $(v_i,v_j)$. By Farkas' lemma the insolvability of $Mx<b$ is equivalent to the existence of a row vector $y = (y_1, ..., y_k)$ satisfying:
\begin{equation}\label{farkas}
y_i\geq0\ \ \forall i,\quad y\neq0,\quad yA=0,\quad y\cdot b\leq0. 
\end{equation}
Exactly the same argument as in the proof of Theorem $14$ of \cite{hopkins} shows that such a $y$ cannot exist. The linear dependencies of $M$ are spanned by sums of rows corresponding to a cycle of $G$. A $y \geq 0 $ satisfying $yA = 0$ would correspond to a sum of potential cycles of $\mathcal{O}$, but each potential cycle has a positive score, so~$y \cdot b > 0$. By construction, $\mathcal{O}_{R} = \mathcal{O}$. 

It remains to be shown that for any region $R$, we have that $\mathcal{O}_R$ is $A$-admissible. Let $R$ be a region of $\Sigma_G(A)$. Encode the inequalities defining $R$ as $Mx < b$. A bad cycle in $\mathcal{O}_R$ corresponds to a vector $y$ satisfying condition (\ref{farkas}) above. Thus, no bad cycles for $\mathcal{O}_R$ can exist. $\square$

\medskip \noindent {\bf Notation.} From now on, for an $A$-admissible partial orientation $\mathcal{O}$, we will use~$\mathrm{rg}(\mathcal{O})$ to denote the unique region of $\Sigma_G(A)$ satisfying $\mathcal{O} = \mathcal{O}_{\mathrm{rg}(\mathcal{O})}$.

\begin{definition}
Let $\mathcal{A}$ be a hyperplane arrangement in $\mathbb{R}^{n}$ and let $W$ be the subspace of $\mathbb{R}^n$ spanned by the normals of the hyperplanes in $\mathcal{A}$. We say that a region $R$ of $\mathcal{A}$ is \emph{relatively bounded} if $R \cap W$ is bounded. The \emph{essentialization} of $\mathcal{A}$ is $\mathcal{A} \cap W$ considered as a hyperplane arrangement in~$W \simeq \mathbb{R}^{k}$, where $k = \mathrm{dim}(W)$.
\end{definition}

The subspace spanned by the normals of the hyperplanes in $\Sigma_G(A)$ is $\mathrm{Span}(\vec{1})^{\perp}$, where $\vec{1} := (1,1,\ldots,1) \in \mathbb{R}^n$.

\begin{thm} \label{thm:bounded}
The relatively bounded regions of $\Sigma_G(A)$ are in bijection with the $A$-admissible partial orientations $\mathcal{O}$ of $G$ for which every step in $\mathcal{O}$ belongs to some potential cycle. The bijection is given by $R \mapsto \mathcal{O}_R$.
\end{thm}

\emph{Proof}: Let $R$ be a region of $\Sigma_G(A)$ and encode the system of inequalities for $R$ as $Mx < b$ as in the proof of Theorem~\ref{thm:farkas}. Recall that the rows of $M$ correspond to steps compatible with $\mathcal{O}_R$. We first need the following:
\begin{quote} 
{\bf Claim}: The region $R$ is relatively bounded if and only if for all vectors~\mbox{$z \lneq 0 \in \mathbb{R}^{k}$}, there exists $y \geq 0 \in \mathbb{R}^{k}$ such that $yM = 0$ and~\mbox{$y \cdot z < 0$}.
\end{quote} 
\noindent \emph{Proof of claim}: By (another version of) Farkas' lemma, the existence of a $y$ as in the claim is equivalent to the nonexistence of a solution $x$ to $Mx \leq z$. So we will show that $R$ is not relatively bounded if and only if there exists~$z \lneq 0$ and $x$ such that $Mx \leq z$. 

The region $R$ is not relatively bounded if and only if there exists~$x \notin \mathrm{Span}( \vec{1} )$ and~\mbox{$x_0 \in \mathbb{R}^n$} such that $M(x_0 + tx) < b$ for all $t \geq 0$, i.e., $tMx < b - Mx_0$ for all~$t \geq 0$. In this case, $Mx \leq 0$. So if $R$ is not relatively bounded, take $x$ as above and let~\mbox{$z := Mx$}. We cannot have $Mx = 0$ since for all $\{v_i,v_j\} \in E$, either~$(v_i,v_j)$ or $(v_j,v_i)$ corresponds to a row of $M$ (or both do), and so $\mathrm{ker}(M) = \mathrm{Span}(\vec{1}) \not \ni x$. 

Conversely, suppose there exists $z \lneq 0$ and $x$ such that $Mx \leq z$. It follows that~\mbox{$x \notin \mathrm{ker}(M) =  \mathrm{Span}(\vec{1})$} and choosing any point $x_0 \in R$, we have
\[Mx_0 < b \Rightarrow 0 < b - Mx_0 \Rightarrow tMx \leq tz < b-Mx_0,\] 
for all $t \geq 0$. Hence, $R$ is not relatively bounded. $\square$

Suppose that $R$ is relatively bounded and let $e \in \mathcal{O}_R$. Define $z \lneq 0 \in \mathbb{R}^{k}$ by 
\[z_{ij} = \begin{cases}
-1& \textrm{ if $e^{-} = v_i$ and $e^{+} = v_j$,}\\
\hfill 0& \textrm{ otherwise.}
\end{cases}\]
Then there must exist $y$ satisfying the condition in the claim. Since, as in the proof of Theorem~\ref{thm:farkas}, the potential cycles of $\mathcal{O}_R$ span the linear dependencies of the rows of $M$, the support of $y$ contains a potential cycle of $\mathcal{O}_R$ containing the step $e$. 

Conversely, let $\mathcal{O}$ be an $A$-admissible partial orientation where each step in~$\mathcal{O}$ is part of a potential cycle. Encode the inequalities of $\mathrm{rg}(\mathcal{O})$ as $Mx < b$. Let~\mbox{$z \lneq 0$} be a vector in $\mathbb{R}^{k}$ and suppose $z_l < 0$ where the $l$th row of $M$ corresponds to a step~\mbox{$e = (v_i,v_j)$}. Let $C$ be a potential cycle containing $e$ (and note that if~\mbox{$e \notin \mathcal{O}$}, then $\{(v_i,v_j),(v_j,v_i)\}$ is a potential cycle for $\mathcal{O}$ containing $e$). The vector corresponding to the steps of $C$ satisfies the condition on $y$ in the claim, and thus $\mathrm{rg}(\mathcal{O})$ is relatively bounded. $\square$

\section{\texorpdfstring{Parking functions and the regions of $\Sigma_G(A)$}{Parking functions and the regions of SG(A)}} \label{sec:pfs}


In this section, we explain how the indegree sequences of partial orientations of~$G$ are closely related to the parking functions of $G_{\sbullet}$, the graph obtained from~$G$ by adding a vertex $v_0$ and an edge between $v_0$ and each~$v \in V$. We will use $V_{\sbullet}$ and~$E_{\sbullet}$ to denote the vertex and edge set, respectively, of $G_{\sbullet}$. Our goal is to show that a natural set of labels for the regions of  $\Sigma_G(A)$ are the set of parking functions of~$G_{\sbullet}$. Let $\mathbb{Z}V$ denote the free abelian group on the vertices in $V$.

\begin{definition}
Let $\mathcal{O}$ be a partial orientation of $G$. The \emph{indegree sequence} of $\mathcal{O}$, denoted $\mathrm{indeg}(\mathcal{O})$, is $\sum_{i=1}^{n} \mathrm{indeg}_{\mathcal{O}}(v_i)\, v_i \in \mathbb{Z}V$.
\end{definition}

\begin{definition}
A \emph{parking function} $c = \sum_{i=1}^{n} c_i v_i$ of $G_{\sbullet}$ with respect to~$v_0$ is an element of $\mathbb{Z}V$ such that for every non-empty subset $W \subseteq V$, there exists $v_i \in W$ with $0 \leq c_i < d_{W}(v_i)$, where $d_{W}(v_i)$ is the number of edges~$e = \{v_i,u\} \in E_{\sbullet}$ with $u \in V_{\sbullet} \setminus W$. For $c, c' \in \mathbb{Z}V$, we write $c \leq c'$ if~$c_i \leq c'_i$ for all $0 \leq i \leq n$. A parking function $c'$ is \emph{maximal} if $c' \leq c$ for any parking function $c$ implies $c = c'$.
\end{definition}

Graphical parking functions were first formally introduced in \cite{postnikov}. However, the essentially equivalent notion of \emph{superstable configurations} has been studied for longer in the context of the abelian sandpile model; see \cite[\S2.4]{hopkins} for a definition of these and their connection to parking functions. One easy observation from the above definition is that if $c'$ is a maximal parking function and $c \in \mathbb{Z}V$ with~\mbox{$0 \leq c \leq c'$}, then $c$ is a parking function as well. The following characterization of maximal parking functions is Theorem~$3.1$ of Benson, Chakrabarty, and Tetali~\cite{benson}:

\begin{thm} A \emph{total orientation} of a graph is a partial orientation where every edge is oriented. A \emph{source} of a total orientation is a vertex whose outdegree equals its degree. The acyclic total orientations of $G_{\sbullet}$ with unique source $v_0$ are in bijection with the maximal parking functions of $G_{\sbullet}$ with respect to $v_0$. The bijection is given by $\mathcal{O} \mapsto \mathrm{indeg}(\mathcal{O}) - \sum_{i=1}^{n} v_i$.
\end{thm}

\begin{prop} \label{prop:parkfns}
The set
\[\{\mathrm{indeg}(\mathcal{O})\colon \mathcal{O} \textrm{ is an acyclic partial orientation of $G$} \}\]
is the set of parking functions of $G_{\sbullet}$ with respect to $v_0$.
\end{prop}

\emph{Proof}: If $\mathcal{O}$ is an acyclic total orientation of $G$, then $\mathcal{O}' := \mathcal{O} \cup \{ (v_0,v_i)\}_{i=1}^n$ is an acyclic total orientation of $G_{\sbullet}$ with unique source $v_0$, and $\mathrm{indeg}(\mathcal{O})$ is equal to $\mathrm{indeg}(\mathcal{O}') - \sum_{i=1}^{n} v_i$. Conversely, any acyclic total orientation~$\mathcal{O}'$ of~$G_{\sbullet}$ with unique source $v_0$ restricts to an acyclic total orientation $\mathcal{O}$ of~$G$, and the same indegree sequence identity holds. So maximal parking functions of $G_{\sbullet}$ correspond to acyclic total orientations of $G$.

Consider any parking function $c$ of $G_{\sbullet}$. There exists a maximal parking function~$c'$ with $c \leq c'$. Let $\mathcal{O}'$ be the acyclic total orientation such that~\mbox{$\mathrm{indeg}(\mathcal{O}') = c'$}. We can easily find $\mathcal{O} \subseteq \mathcal{O}'$ so that $\mathrm{indeg}(\mathcal{O}) = c$. Conversely, consider some acyclic partial orientation $\mathcal{O}$ of $G$.  Take any edge~$e = \{v_i,v_j\}$ that is blank in~$\mathcal{O}$. Suppose both $\mathcal{O} \cup \{(v_i,v_j)\}$ and $\mathcal{O} \cup \{(v_j,v_i)\}$ have directed cycles. Then there was a directed path from $v_i$ to $v_j$ in~$\mathcal{O}$ and a directed path from $v_j$ to $v_i$ in $\mathcal{O}$, so $\mathcal{O}$ already had a directed cycle, yielding a contradiction. Thus, we can repeatedly orient blank edges in $\mathcal{O}$ and arrive at an acyclic total orientation $\mathcal{O}' \supseteq \mathcal{O}$. Then~\mbox{$0 \leq \mathrm{indeg}(\mathcal{O}) \leq \mathrm{indeg}(\mathcal{O}')$} and $\mathrm{indeg}(\mathcal{O}')$ is a maximal parking function of~$G_{\sbullet}$. Therefore, $\mathrm{indeg}(\mathcal{O})$ is indeed a parking function. $\square$

Theorem~\ref{thm:farkas} makes clear that an $A$-admissible partial orientation is necessarily acyclic, as follows. Suppose $\mathcal{O}$ is a partial orientation that contains a directed cycle of steps $C = \{(u_1,u_2),(u_2,u_3),\ldots,(u_k,u_1)\}$. Then let
\[C^{\,\mathrm{rev}}:=\{(u_1,u_k),(u_k,u_{k-1}),\ldots,(u_2,u_1)\}.\] 
Since $\Sigma_G(A)$ has a central region $R_0$, the partial orientation $\mathcal{O}_{R_0} = \emptyset$ is $A$-admissible. So we have~$\nu(C,\mathcal{O}) = -\nu(C^{\,\mathrm{rev}},\mathcal{O}_{R_0}) < 0$ and thus $\mathcal{O}$ is not $A$-admissible. Thus, if~$\mathcal{O}$ is $A$-admissible, Proposition~\ref{prop:parkfns} implies that $\mathrm{indeg}(\mathcal{O})$ is a parking function of $G_{\sbullet}$. What remains to be proven is that any parking function of $G_{\sbullet}$ can be realized as $\mathrm{indeg}(\mathcal{O})$ for some $A$-admissible $\mathcal{O}$. We are therefore interested in building $A$-admissible partial orientations with particular indegree sequences. The following ``topological'' lemma will allow us to build up $A$-admissible partial orientations from other $A$-admissible partial orientations with some control over the resulting indegrees.

\begin{lemma} \label{lem:top}
Let $\mathcal{O}$ be an $A$-admissible partial orientation, and let $W \subseteq V$ be a subset of the vertices of $G$ satisfying:
\begin{enumerate}
\item there do not exist $u \in W^{c}$, $w \in W$ with $(w,u) \in \mathcal{O}$;\label{cond:noarrowsout}
\item there is some $u \in W^{c}$, $w \in W$ such that $\{u,w\}$ is a blank edge of $\mathcal{O}$.\label{cond:ablank}
\end{enumerate}
Then there exists $u \in W^{c}$ and $w \in W$ such that $\mathcal{O} \, \cup \{(u,w)\}$ is also $A$-admissible.
\end{lemma}

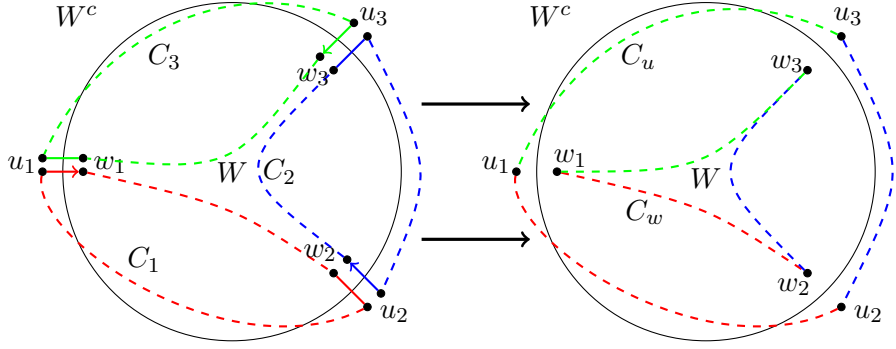
\begin{figure}
\begin{center}
\begin{tikzpicture}[scale=0.9]

\draw(0,0) circle(2.5);
\node at (0,0) {$W$};
\node at (-2.3,2.3) {$W^c$};
\node (u1) at (-2.8,0) [circle,draw,fill=black,inner sep=1pt] {};
\node (w1) at (-2.2,0) [circle,draw,fill=black,inner sep=1pt] {};
\draw [->,color=red,style=thick] (u1) -- (w1);
\node (u2) at (2,-2) [circle,draw,fill=black,inner sep=1pt] {};
\node (w2) at (1.5,-1.5) [circle,draw,fill=black,inner sep=1pt] {};
\draw [color=red,style=thick] (w2) -- (u2);

\node at (-3.1,0.1) {$u_1$};
\node at (-1.8,0.1) {$w_1$};

\draw [color=red,style=dashed,style=thick] (w1) .. controls (0,-0.5) .. (w2);
\draw [color=red,style=dashed,style=thick] (u2) .. controls (0,-3) and (-3,-1) .. (u1);

\node (u2p) at (2.2,-1.8) [circle,draw,fill=black,inner sep=1pt] {};
\node (w2p) at (1.7,-1.3) [circle,draw,fill=black,inner sep=1pt] {};
\draw [->,color=blue,style=thick] (u2p) -- (w2p);
\node (u3) at (2,2) [circle,draw,fill=black,inner sep=1pt] {};
\node (w3) at (1.5,1.5) [circle,draw,fill=black,inner sep=1pt] {};
\draw [color=blue,style=thick] (w3) -- (u3);

\node at (2.4,-2.1) {$u_2$};
\node at (1.3,-1.2) {$w_2$};

\draw [color=blue,style=dashed,style=thick] (w2p) .. controls (0,0) .. (w3);
\draw [color=blue,style=dashed,style=thick] (u3) .. controls (3,0) .. (u2p);

\node (u3p) at (1.8,2.2) [circle,draw,fill=black,inner sep=1pt] {};
\node (w3p) at (1.3,1.7) [circle,draw,fill=black,inner sep=1pt] {};
\draw [->,color=green,style=thick] (u3p) -- (w3p);
\node (u1p) at (-2.8,0.2) [circle,draw,fill=black,inner sep=1pt] {};
\node (w1p) at (-2.2,0.2) [circle,draw,fill=black,inner sep=1pt] {};
\draw [color=green,style=thick] (w1p) -- (u1p);

\node at (2.1,2.3) {$u_3$};
\node at (1.2,1.4) {$w_3$};

\draw [color=green,style=dashed,style=thick] (w3p) .. controls (0,0) .. (w1p);
\draw [color=green,style=dashed,style=thick] (u3p) .. controls (0,3) and (-2,2) .. (u1p);

\node at (-1.3,-1.3) {$C_1$};
\node at (0.7,0) {$C_2$};
\node at (-1,1.7) {$C_3$};

\draw [->,style=very thick] (2.8,1) -- (4.4,1);
\draw [->,style=very thick] (2.8,-1) -- (4.4,-1);

\draw(7,0) circle(2.5);
\node at (7,-0.1) {$W$};
\node at (7-2.3,2.3) {$W^c$};
\node (u1) at (7+-2.8,0) [circle,draw,fill=black,inner sep=1pt] {};
\node (w1) at (7+-2.2,0) [circle,draw,fill=black,inner sep=1pt] {};
\node (u2) at (7+2,-2) [circle,draw,fill=black,inner sep=1pt] {};
\node (w2) at (7+1.5,-1.5) [circle,draw,fill=black,inner sep=1pt] {};
\node (u3) at (7+2,2) [circle,draw,fill=black,inner sep=1pt] {};
\node (w3) at (7+1.5,1.5) [circle,draw,fill=black,inner sep=1pt] {};

\node at (7+-3.1,0.1) {$u_1$};
\node at (7+-2.0,0.2) {$w_1$};

\node at (7+2.4,-2.1) {$u_2$};
\node at (7+1.3,-1.7) {$w_2$};

\node at (7+2.1,2.3) {$u_3$};
\node at (7+1.2,1.6) {$w_3$};

\draw [color=red,style=dashed,style=thick] (w1) .. controls (7+0,-0.5) .. (w2);
\draw [color=red,style=dashed,style=thick] (u2) .. controls (7+0,-3) and (7+-3,-1) .. (u1);

\draw [color=blue,style=dashed,style=thick] (w2) .. controls (7+0,0) .. (w3);
\draw [color=blue,style=dashed,style=thick] (u3) .. controls (7+3,0) .. (u2);

\draw [color=green,style=dashed,style=thick] (w3) .. controls (7+0,0) .. (w1);
\draw [color=green,style=dashed,style=thick] (u3) .. controls (7+0,3) and (7+-2,2) .. (u1);

\node at (7-0.9,-0.6) {$C_w$};
\node at (7+-1,1.7) {$C_u$};

\end{tikzpicture}
\end{center}
\caption{Diagram explaining the equality of the score sums in the proof of Lemma~\ref{lem:top} with $i = 1$ and $j = 4$.} \label{fig:top}

\end{figure}

\emph{Proof}: In this proof, all $u_k$ are elements of $W^{c}$ and all $w_k$ are elements of $W$. By (\ref{cond:ablank}), we may choose a blank edge $\{u_1,w_1\}$ of $\mathcal{O}$. If $\mathcal{O}_1 :=\mathcal{O}\cup\{(u_1,w_1)\}$ is not $A$-admissible, there exists some bad potential cycle $C_1$ containing $(u_1,w_1)$ and some $(w_2,u_2)$, where $\{u_2,w_2\}$ is necessarily blank by (\ref{cond:noarrowsout}). Next we consider~\mbox{$\mathcal{O}_2 := \mathcal{O} \cup \{(u_2,w_2)\}$}; if this not admissible, we get a bad cycle $C_2$ containing~$(u_2,w_2)$ and $(w_3,u_3)$, and so on. Either we arrive at an admissible partial orientation, or this process goes on forever. Suppose it goes on forever. Because there are only a finite number of blank edges between $W^{c}$ and $W$, eventually we obtain $i < j$ where $(u_i,w_i) = (u_j,w_j)$. Consider $\sum_{k=i}^{j-1} \nu(C_k, \mathcal{O}_k)$. In this sum, the contribution of the step~$(u_k,w_k) \in \mathcal{O}_k$ in $C_k$ cancels with the contribution of the step~$(w_k,u_k)~\notin~\mathcal{O}_{k-1}$ in $C_{k-1}$, and what remains is the sum of the scores relative to $\mathcal{O}$ of a cycle $C_u$ joining $u_{j-1}$ to $u_{j-2}$ and so on to $u_{i}$ and back to~$u_{j-1}$ and of a cycle $C_w$ joining $w_i$ to $w_{i+1}$ and so on to $w_{j-1}$ and back to $w_i$. Neither of the cycles $C_u$ or $C_w$ contain any of the directed edges~$(u_k,w_k)$. Figure~\ref{fig:top} gives a diagrammatic explanation of the equality of these cycle score sums. But then~\mbox{$0 < \nu(C_u,\mathcal{O}) + \nu(C_w,\mathcal{O}) = \sum_{k=i}^{j-1} \nu(C_k, \mathcal{O}_k) \leq 0$}, a contradiction. $\square$

\begin{thm} \label{thm:admissibles}
Let $\mathcal{O}$ be an acyclic partial orientation of $G$. Then there exists an $A$-admissible partial orientation $\mathcal{O}'$ such that $\mathrm{indeg}(\mathcal{O}) = \mathrm{indeg}(\mathcal{O}')$.
\end{thm}

\emph{Proof}: If $\mathcal{O} = \emptyset$, there is nothing to prove. Otherwise, set $\mathcal{O}_0 := \emptyset$ and recursively define $\mathcal{O}_i$ from $\mathcal{O}_{i-1}$ as follows. Let $W_i$ be the set of $v \in V$ such that~\mbox{$\mathrm{indeg}_{\mathcal{O}_{i-1}}(v) < \mathrm{indeg}_{\mathcal{O}}(v)$}. If $W_i \neq \emptyset$, apply Lemma~\ref{lem:top} to $W_i$ and $\mathcal{O}_{i-1}$ and let $\mathcal{O}_{i}$ be the resulting partial orientation, i.e., set $\mathcal{O}_i := \mathcal{O}_{i-1} \cup \{(u,w)\}$ for the appropriate $u \in W_i^{c}$ and~$w \in W_i$. We claim that the lemma may indeed be applied at every iteration. Suppose not. The first condition of the lemma clearly applies by construction. So suppose there exists $1 \leq i \leq |\mathcal{O}|$ such that every edge between some $u \in W_i^{c}$ and some $w \in W_i$ is already oriented as $(u,w)$ in $\mathcal{O}_{i-1}$. Then for each~$w \in W_i$, we have $\mathrm{indeg}_{\mathcal{O}}(w) > |\{ \{u,w\} \in E: u \in W_i^C\}|$, so $w$ must have an arrow in $\mathcal{O}$ coming into it from some other vertex in $W_i$. But this forces $\mathcal{O}$ to contain a cycle of steps involving the vertices of $W_i$, which is a contradiction. Therefore Lemma~\ref{lem:top} applies at every iteration as claimed. Setting $\mathcal{O}' := \mathcal{O}_{|\mathcal{O}|}$, we arrive at an $A$-admissible partial orientation with the desired indegree sequence. $\square$ 

We are now prepared to prove the main result of this section, Corollary~\ref{cor:main}, which establishes a conjecture of Hopkins and Perkinson~\cite{hopkins}. Indeed, Corollary~\ref{cor:main} subsumes the main result of that paper (which was proved in a different way using the abelian sandpile model).

\begin{definition} The following procedure is called the \emph{Pak-Stanley labeling} of a bigraphical arrangement. It labels each region with an element of $\mathbb{Z}V$. Label the central region of $\Sigma_G(A)$ with $0$.  Put the central region in a queue,~$Q$.  Then, as long as $Q$ is not empty:
\begin{enumerate}
  \item Remove the first region $R$ from $Q$.
  \item For each unlabeled region $R'$ bordering~$R$:
    \begin{enumerate}
      \item Determine the unique indices $i\neq j$ such that $x_j-x_i<a_{ji}$ in~$R$ but~$x_j > x_i + a_{ji}$ in $R'$.
  \item If $R$ is labeled by  $c=\sum_{k=1}^nc_kv_k$, then label $R'$ by $c'=c+v_j$.  
  \item Add $R'$ to the end of $Q$.
  \end{enumerate}
\end{enumerate}
Let $\lambda(R)$ denote the Pak-Stanley label of the region $R$.
\end{definition}

\begin{cor} \label{cor:main}
The set
\[\{\lambda(R)\colon \textrm{$R$ is a region of $\Sigma_G(A)$}\}\]
is the set of all parking functions of $G_{\sbullet}$ with respect to $v_0$.
\end{cor}

\emph{Proof}: We first check inductively that $\lambda(R) = \mathrm{indeg}(\mathcal{O}_R)$. This identity clearly holds for the central region. So suppose $\lambda(R) = \mathrm{indeg}(\mathcal{O}_R)$ and that~$R'$ borders~$R$ with $x_j-x_i<a_{ji}$ in $R$ but $x_j > x_i + a_{ji}$ in $R'$. Then we have
\[\lambda(R') = \lambda(R) + v_j = \mathrm{indeg}(\mathcal{O}_R) + v_j = \mathrm{indeg}(\mathcal{O}_{R'}).\]

By Proposition~\ref{prop:parkfns}, for any parking function $c$ of $G_{\sbullet}$, there is some acyclic partial orientation $\mathcal{O}$ such that $\mathrm{indeg}(\mathcal{O}) = c$. Finally, by Theorem~\ref{thm:admissibles}, we can find an $A$-admissible orientation $\mathcal{O}'$ with the same indegree sequence as~$\mathcal{O}$, so we have~$\lambda(\mathrm{rg}(\mathcal{O}'))  = c$. $\square$

\begin{remark}\label{remark:parkfns}
In~\cite[\S4]{hopkins}, it is shown that
\begin{align*}
\{ c \in \mathbb{Z}V \colon &\textrm{$c = \mathrm{indeg}(\mathcal{O}) -  \sum_{i=1}^{n} v_i$ for some acyclic partial orientation $\mathcal{O}$ of $G$,} \\
&\textrm{$c_i = -1$, and $c_j \geq 0$ for all $j \neq i$} \}
\end{align*}
is the set of parking functions of $G$ with respect to $v_i$. 
Thus, in light of Theorem~\ref{farkas}, Proposition~\ref{prop:parkfns}, and
Theorem~\ref{thm:admissibles}, the parking functions of $G$ with respect to each of
its vertices are encoded in the regions of $\Sigma_G(A)$.
\end{remark} 

\begin{cor} \label{cor:minreg}
The number of regions of $\Sigma_G(A)$ is at least the number of spanning trees of $G_{\sbullet}$.
\end{cor}

\emph{Proof}: The number of spanning trees of a graph equals the number of parking functions of that graph with respect to any vertex (see~\cite{chebikin} or \cite{postnikov}), so Corollary~\ref{cor:main} implies this lower bound. $\square$

\section{\texorpdfstring{Number of regions of $\Sigma_G(A)$}{Number of regions of SG(A)}} \label{sec:numregs}

We have already seen that the number of regions of $\Sigma_G(A)$ is at least the number of spanning trees of $G_{\sbullet}$. In this section we give further bounds on the number of regions of $\Sigma_G(A)$. The graph $G$ remains fixed, but we will allow the parameter list~$A$ to vary (while always maintaining a central region). 

\medskip

\noindent {\bf Notation.} We will denote the number of regions of $\Sigma_G(A)$ by $r(\Sigma_G(A))$ and the number of relatively bounded regions by $b(\Sigma_G(A))$.

\begin{definition}
Let $\mathcal{A}$ be a hyperplane arrangement. The hyperplanes $H_{1}, \ldots, H_{k}$ in $\mathcal{A}$ are \emph{linearly independent} if their normals are linearly independent. The arrangement $\mathcal{A}$ is \emph{generic} if
\[ H_{1} \cap \cdots \cap H_{k} \neq \emptyset \quad \Leftrightarrow \quad H_{1}, \ldots, H_{k} \text{ are linearly independent} \]
for all subsets $\{H_{1}, \ldots, H_{k}\} \subseteq \mathcal{A}$.
\end{definition}

For instance, $\Sigma_G(A)$ is generic when the $a_{ij}$ are algebraically independent. We will use $\mathrm{GEN}$ to denote the parameter list of an arbitrary generic bigraphical arrangement.

\begin{thm} \label{thm:genreg}
For any generic bigraphical arrangement, the characteristic polynomial of $\Sigma_G(\mathrm{GEN})$ is given by
\[\bigchi_{\Sigma_G(\mathrm{GEN})}(t) = (-2)^{n-1} \; t \; T_G(1-t/2,1),\]
where $T_G(x,y)$ is the Tutte polynomial of $G$. Consequently,
\begin{align*}
r(\Sigma_G(\mathrm{GEN})) &= 2^{n-1}T_G(3/2,1); \\
b(\Sigma_G(\mathrm{GEN})) &= 2^{n-1}T_G(1/2,1).
\end{align*}
\end{thm}

\emph{Proof}: As Stanley shows in~\cite[p.~412]{stanley}, the characteristic polynomial of a generic arrangement $\mathcal{A}$ is given by
\[\bigchi_{\mathcal{A}}(t) = \sum_{\mathcal{B}} (-1)^{|\mathcal{B}|} t^{n-|\mathcal{B}|},\]
where the sum is over all linearly independent subsets $\mathcal{B}$ of $\mathcal{A}$. So we must find the linearly independent subsets of $\Sigma_G(\mathrm{GEN})$. Let the hyperplanes of~$\Sigma_G(\mathrm{GEN})$ be $H_e$ with linear parts $L_e := x_i - x_j$ for steps $e = (v_i, v_j)$ of~$G$. For a step $e = (v_i,v_j)$ of~$G$, define $\pi(e) = \{v_i,v_j\}$. For $\mathcal{B} \subseteq \Sigma_G(\mathrm{GEN})$, let $\pi(\mathcal{B})$ be the multiset $\{\pi(e)\}_{H_e \in \mathcal{B}}$. We claim that there exists a linear dependence among $\mathcal{B}$ if and only if there exists a cycle of undirected edges $\{ \{u_1,u_2\}, \{u_2, u_3\}, \ldots, \{u_k, u_1\} \} \subseteq \pi(\mathcal{B})$. Suppose there exists such a cycle and $H_{e_1}, \ldots, H_{e_k}$ are the corresponding hyperplanes. Then set~$\lambda_{e_i} = 1$ if~$e_i = (u_i, u_j)$ and $\lambda_{e_i} = -1$ if $e_i = (u_j,u_i)$. We see that~\mbox{$\sum_{i=1}^{k} \lambda_{e_i} L_{e_i} = 0$}. Conversely, suppose there exists a linear dependence in some subset~$\mathcal{B} \subseteq \Sigma_G(\mathrm{GEN})$. Similarly to the proof of Theorem~\ref{thm:farkas}, the undirected cycles of $G$ span the linear dependencies of $\Sigma_G(\mathrm{GEN})$, so $\mathcal{B}$ must contain such a cycle. Thus,~$\mathcal{B}$ is linearly independent if and only if $\pi(\mathcal{B})$ is a forest of $G$. For each edge~$e = \{v_i,v_j\}$ in such forest, there are two hyperplanes $H_{(v_i,v_j)}$ or $H_{(v_j,v_i)}$ we could include in $\mathcal{B}$. A forest~$F$ thus corresponds to $2^{|F|}$ linearly independent subsets $\mathcal{B} \subseteq \Sigma_G(\mathrm{GEN})$. Therefore, the characteristic polynomial of~$\Sigma_G(\mathrm{GEN})$ is
\[\bigchi_{\Sigma_G(\mathrm{GEN})}(t) = \sum_{F} (-2)^{|F|} t^{n-|F|},\]
where the sum is over all forests $F$ of $G$. The following formula (see \cite[p.~1135]{welsh}) relates the Tutte polynomial to the forests of $G$:
\[ \sum_{i=0}^{n-1} f_i \, t^i = t^{n-1}T_G(1 + 1/t,1),\]
where $f_i$ is the number of forests of $G$ of size $i$. So we have
\begin{align*}
 \bigchi_{\Sigma_G(\mathrm{GEN})}(t) &= \sum_{i=0}^{n-1} f_i (-2)^{i} t^{n-i} \\
 &= t^n \sum_{i=0}^{n-1} f_i \left(\frac{-2}{t}\right)^{i} \\
 &= (-2)^{n-1} \; t \; T_G(1-t/2,1).
\end{align*}
A classical result in the theory of hyperplane arrangements, Zaslavsky's theorem~\cite{zaslavsky}, relates the number of regions and number of relatively bounded regions of a hyperplane arrangement $\mathcal{A}$ to its characteristic polynomial:
\begin{align*}
r(\mathcal{A}) &= |\bigchi_{\mathcal{A}}(-1)|; \\
b(\mathcal{A}) &= |\bigchi_{\mathcal{A}}(1)|.
\end{align*}
Zavlasky's theorem establishes that the formulas for the number of regions and bounded regions are as claimed. $\square$

\begin{cor} \label{cor:joke}
Suppose $G$ is planar and $G^{*}$ is its dual graph. Then the following are equal:
\begin{enumerate}
\item the probability that after removing each edge from $G^{*}$ with probability~$2/3$, the resulting graph remains connected;
\item the probability that a partial orientation of $G$ chosen uniformly at random is $\mathrm{GEN}$-admissible.
\end{enumerate}
\end{cor}

\emph{Proof}: The first probability is given by $R_{G^{*}}(2/3)$, where $R_{G^{*}}(p)$ is the all-terminal reliability polynomial of $G^{*}$. Let $V^{*}$ be the vertex set of $G^{*}$, let~$E^{*}$ be its edge set, and let $k(G^{*})$ be the number of connected components of~$G^{*}$. A well-known formula connecting the reliability and Tutte polynomials (see~\cite[p.~1335]{welsh}) is 
\[R_{G^{*}}(p) = (1-p)^{|V^{*}| - k(G^{*})} p^{|E^{*}|-|V^{*}|+k(G^{*})}T_{G^*}(1,1/p).\] 
Suppose $G$ has $f$ faces. Then the probability of the first event is
  \begin{align*}
  R_{G^{*}}(2/3) &= (1-2/3)^{|V^{*}| - k(G^{*})} (2/3)^{|E^{*}|-|V^{*}|+k(G^{*})}T_{G^*}(1,3/2) \\[5pt]
  &= \frac{2^{|E|-f+1}T_G(3/2,1)}{3^{|E|}} \\[5pt]
  &= \frac{2^{n-1}T_G(3/2,1)}{3^{|E|}} \\[5pt]
  &= \frac{r(\Sigma_{G}(\mathrm{GEN}))}{3^{|E|}},
  \end{align*}
  where in the second line we have used the formula $T_{G^{*}}(y,x) = T_G(x,y)$ (see~\cite[p.1131]{welsh}), in the third line we have used Euler's formula $n - |E| + f =~2$, and in the last line we applied Theorem~\ref{thm:genreg}. The result now follows from Theorem~\ref{thm:farkas}. $\square$
 

\begin{definition}
Let $\mathcal{O}$ be a partial orientation. A potential cycle for $\mathcal{O}$ is \emph{very bad} if has a negative score. We say $\mathcal{O}$ is \emph{almost-$A$-admissible} if it is not $A$-admissible but it has no very bad potential cycles. We say $\mathcal{O}$ is \emph{far-from-$A$-admissible} if it has a very bad potential cycle.
\end{definition}

\begin{prop} \label{prop:zeroreg}
For an almost-$A$-admissible partial orientation $\mathcal{O}$, define~$w(\mathcal{O})$ to be the number of steps $e$ of $G$, including blanks of $\mathcal{O}$, belonging to some potential cycle $C$ for $\mathcal{O}$ with $\nu_{A}(C,\mathcal{O}) = 0$. Define $z(\mathcal{O})$ to be the maximum $k$ such that there exists disjoint potential cycles $C_1,\ldots,C_k$ for $\mathcal{O}$ with $\nu_A(C_i,\mathcal{O}) = 0$ for all $i$. Then
\[ \sum_{\mathcal{O}} \frac{1}{2^{w(\mathcal{O})}} \leq r(\Sigma_G(\mathrm{GEN})) -  r(\Sigma_G(A)) \leq \sum_{\mathcal{O}} \frac{1}{2^{z(\mathcal{O})}}, \]
where the sum is over all almost-$A$-admissible partial orientations $\mathcal{O}$.
\end{prop}
\emph{Proof}: Let $S$ denote the set of $2|E|$ steps of $G$.  For each $(v_i,v_j)\in S$, let~$\varepsilon_{ij} > 0$ be a real number.  For each sign pattern $\sigma \in\{-1,1\}^S$, let $\sigma_{ij}:=\sigma(v_i,v_j)$, and define the parameter list $A^{\sigma}$ with parameters $a_{ij}^{\sigma}:=a_{ij}+{\sigma}_{ij}\,\varepsilon_{ij}$ for all steps~$(v_i,v_j)\in S$.  Take the real numbers~$\varepsilon_{ij}$ small enough and generic so that for all choice of $\sigma$, each $\Sigma_G(A^{\sigma})$ is generic, each $A$-admissible partial orientation is $A^{\sigma}$-admissible, and each far-from-$A$-admissible partial orientation is far-from-$A^{\sigma}$-admissible. No partial orientations are almost-admissible for a generic arrangement. Let $\mathrm{Gain}(\sigma)$ be the set of almost-$A$-admissible partial orientations that are $A^{\sigma}$-admissible. Note that $|\mathrm{Gain}(\sigma)| = r(\Sigma_G(A^{\sigma})) - r(\Sigma_G(A))$. But then by Theorem~\ref{thm:genreg}, we have $r(\Sigma_G(A^{\sigma})) =  r(\Sigma_G(\mathrm{GEN}))$ for any~$\sigma \in \{-1,1\}^S$. So in order to bound $r(\Sigma_G(\mathrm{GEN})) - r(\Sigma_G(A))$, we bound $|\mathrm{Gain}(\sigma)|$.

Define $X = \{ (\sigma,\mathcal{O}): \mathcal{O} \in \mathrm{Gain}(\sigma)\}$. Because $|\mathrm{Gain}(\sigma)|$ is independent of $\sigma$, we have $|X| = 2^{|S|}|\mathrm{Gain}(\sigma)|$. Fix some almost-$A$-admissible partial orientation $\mathcal{O}$ and suppose $w(\mathcal{O}) = k$. Let $S' \subseteq S$ be the set of steps of $G$ belonging to some potential cycle $C$ for $\mathcal{O}$ with $\nu_{A}(C,\mathcal{O}) = 0$. Suppose $\sigma$ is such that for each $e = (v_i,v_j) \in S'$, we have $\sigma_{ij} = 1$ if $e \notin \mathcal{O}$ and $\sigma_{ji} = -1$ if $e \in \mathcal{O}$. Then $\mathcal{O}$ is $A^{\sigma}$ admissible: let $C$ be any potential cycle for $\mathcal{O}$ with~$\nu_{A}(C,\mathcal{O}) = 0$; then,
\begin{align*}
\nu_{A^{\sigma}}(C,\mathcal{O}) &= \nu_{A}(C,\mathcal{O}) + \hspace{-0.3cm} \sum_{\substack{(v_i,v_j) \in C \\ \{v_i,v_j\}  \notin \mathcal{O}}}  \hspace{-0.3cm} \sigma_{ij}\varepsilon_{ij} -  \hspace{-0.3cm} \sum_{\substack{(v_i,v_j) \in C \\(v_i,v_j) \in \mathcal{O}}}  \hspace{-0.3cm} \sigma_{ji}\varepsilon_{ji}\\
&=  \sum_{\substack{(v_i,v_j) \in C \\ \{v_i,v_j\}  \notin \mathcal{O}}}  \hspace{-0.3cm}\varepsilon_{ij} +  \hspace{-0.3cm} \sum_{\substack{(v_i,v_j) \in C \\(v_i,v_j) \in \mathcal{O}}}  \hspace{-0.3cm} \varepsilon_{ji} > 0.
\end{align*}
Since we are free to choose the sign of $\sigma$ associated to any step not in~$S'$, there are at least $2^{|S|-k}$ sign patterns $\sigma$ with $\mathcal{O} \in \mathrm{Gain}(\sigma)$. We have $|\{\sigma: (\sigma,\mathcal{O}) \in X\}| \geq 2^{|S|-k}.$
So, 
\[
|X| = \sum_{\mathcal{O}} |\{\sigma: (\sigma,\mathcal{O}) \in X\}| \geq \sum_{\mathcal{O}} 2^{|S|-w(\mathcal{O})},
\]
where the sum is over all almost-$A$-admissible partial orientations $\mathcal{O}$.

Now fix some almost-$A$-admissible partial orientation $\mathcal{O}$ and suppose $z(\mathcal{O}) = k$. Let $\mathcal{C} := \{C_1,\ldots,C_k\}$ be a set of disjoint potential cycles for~$\mathcal{O}$ with $\nu_{A}(C_i,\mathcal{O}) = 0$ for all $i$. For each $\sigma \in \{-1,1\}^{S}$ and each subset~$\mathcal{D} \subseteq \mathcal{C}$, define $\mathrm{flip}(\sigma,\mathcal{D})$ to be the same sign pattern as $\sigma$ except that for any step~$e = (v_i,v_j)$ in a cycle $C \in \mathcal{D}$, we have $\mathrm{flip}(\sigma,\mathcal{D})_{ij} = -\sigma_{ij}$ if $e \notin \mathcal{O}$ and $\mathrm{flip}(\sigma,\mathcal{D})_{ji} = -\sigma_{ji}$ if $e \in \mathcal{O}$. Write $\sigma \sim \tau$ if there exits $\mathcal{D} \subseteq \mathcal{C}$ such that $\tau = \mathrm{flip}(\sigma,\mathcal{D})$. The relation $\sim$ defines an equivalence relation. Further,~$\mathcal{O}$ belongs to at most one of the $\mathrm{Gain}(\sigma)$ among all $\sigma$ in some equivalence class. Indeed, suppose $\mathcal{O} \in \mathrm{Gain}(\sigma)$ and $\tau \sim \sigma$ but $\tau \neq \sigma$. Then let~$\mathcal{D}$ be such that $\tau = \mathrm{flip}(\sigma,\mathcal{D})$. Note that $\mathcal{D}$ is nonempty. So for any potential cycle~$C \in \mathcal{D}$, we have
\begin{align*}
\nu_{A^{\tau}}(C,\mathcal{O}) &= \nu_{A}(C,\mathcal{O}) + \hspace{-0.3cm} \sum_{\substack{(v_i,v_j) \in C \\ \{v_i,v_j\}  \notin \mathcal{O}}}  \hspace{-0.3cm} \tau_{ij}\varepsilon_{ij} -  \hspace{-0.3cm} \sum_{\substack{(v_i,v_j) \in C \\(v_i,v_j) \in \mathcal{O}}}  \hspace{-0.3cm} \tau_{ji}\varepsilon_{ji}\\
&= - \hspace{-0.3cm} \sum_{\substack{(v_i,v_j) \in C \\ \{v_i,v_j\} \notin \mathcal{O}}} \hspace{-0.3cm} \sigma_{ij}\varepsilon_{ij} + \hspace{-0.3cm} \sum_{\substack{(v_i,v_j) \in C \\(v_i,v_j) \in \mathcal{O}}} \hspace{-0.3cm} \sigma_{ji}\varepsilon_{ji}\\
&= -\nu_{A^{\sigma}}(C,\mathcal{O}) < 0.
\end{align*}
Thus, $\mathcal{O} \notin \mathrm{Gain}(\tau)$. Since $\mathcal{O}$ is in at most one of the $\mathrm{Gain}(\sigma)$ among~$\sigma$ in an equivalence class and each class has $2^k$ members, we have \mbox{$|\{\sigma: (\sigma,\mathcal{O}) \in X\}| \leq 2^{|S|-k}$}.
So, 
\[
|X| = \sum_{\mathcal{O}} |\{\sigma: (\sigma,\mathcal{O}) \in X\}| \leq \sum_{\mathcal{O}} 2^{|S|-z(\mathcal{O})},
\]
where the sum is over all almost-$A$-admissible partial orientations $\mathcal{O}$.  But then recall that~\mbox{$|\mathrm{Gain}(\sigma)| = (1/2^{|S|})|X|$}. Therefore, 
\[  \sum_{\mathcal{O}} \frac{1}{2^{w(\mathcal{O})}} \leq  |\mathrm{Gain}(\sigma)| \leq \sum_{\mathcal{O}} \frac{1}{2^{z(\mathcal{O})}}. \; \square\]

\begin{example} \label{ex:semicyc}
Consider the cycle graph $C_n$, labeled as below:
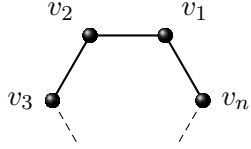
\begin{center}
\begin{tikzpicture}
\SetVertexMath
\GraphInit[vstyle=Art]
\tikzset{VertexStyle/.style = {%
shape = circle,
shading = ball,
ball color = black,
inner sep = 2pt
}}
\SetUpVertex[MinSize=3pt]
\SetVertexLabel
\SetUpEdge[color=black]
\Vertex[LabelOut,Lpos=0,Ldist=.01cm,x=1,y=0]{v_n}
\Vertex[LabelOut,Lpos=45,Ldist=.01cm,x=0.5,y=0.866]{v_1}
\Vertex[LabelOut,Lpos=135,Ldist=.01cm,x=-0.5,y=0.866]{v_2}
\Vertex[LabelOut,Lpos=180,Ldist=.01cm,x=-1,y=0]{v_3}
\Edge(v_n)(v_1)
\Edge(v_1)(v_2)
\Edge(v_2)(v_3)
\draw[style=densely dashed] (-1,0) -- (-0.65,-0.606);
\draw[style=densely dashed] (1,0) -- (0.65,-0.606);
\end{tikzpicture}
\end{center}
\smallskip
The right inequality from Proposition~\ref{prop:zeroreg} becomes an equality for $C_n$: for any almost-$A$-admissible $\mathcal{O}$, there is only one potential cycle $C$ for $\mathcal{O}$ for which we have~\mbox{$\nu_{A}(C,\mathcal{O}) = 0$}. Thus, using the notation from the proof of Proposition~\ref{prop:zeroreg}, we have $\mathcal{O} \in \mathrm{Gain}(\sigma)$ if and only if $\mathcal{O} \notin \mathrm{Gain}(\mathrm{flip}(\sigma,\{C\}))$. Further, any $\mathcal{O} \in \mathrm{Gain}(\sigma)$ is relatively bounded by the condition of Theorem~\ref{thm:bounded}. So,
\begin{align*}
 r(\Sigma_{C_n}(A)) &= r(\Sigma_{C_n}(\mathrm{GEN})) - |\{\mathcal{O}\colon \mathcal{O} \textrm{ is almost-$A$-admissible}\}|/2; \\
 b(\Sigma_{C_n}(A)) &= b(\Sigma_{C_n}(\mathrm{GEN})) - |\{\mathcal{O}\colon \mathcal{O} \textrm{ is almost-$A$-admissible}\}|/2.
\end{align*}
This lets us compute exact formulas for $r(\Sigma_{C_n}(A))$. For instance,
\begin{align*}
 r(\Sigma_{C_n}(\mathrm{SEMI})) &= \begin{cases} 3^n - 2^n &\textrm{ if $n$ is odd,} \\ 3^n - 2^n - \binom{n}{n/2} &\textrm{ if $n$ is even.} \end{cases} \\
 b(\Sigma_{C_n}(\mathrm{SEMI})) &= \begin{cases} 2^n - 1 &\textrm{ if $n$ is odd,} \\ 2^n - 1 - \binom{n}{n/2} &\textrm{ if $n$ is even.} \end{cases} \\
 r(\Sigma_{C_n}(\mathrm{SHI})) &= 3^n - 2^n - n. \\
 b(\Sigma_{C_n}(\mathrm{SHI})) &= 2^n - 1 - n.
\end{align*}
To see this, first note that $T_{C_n}(x,y) = y + \sum_{i=1}^{n-1} x^i$, so
\[r(\Sigma_{C_n}(\mathrm{GEN})) = 2^{n-1}T_{C_n}(3/2,1) = 3^n - 2^n;\]
\[b(\Sigma_{C_n}(\mathrm{GEN})) = 2^{n-1}T_{C_n}(1/2,1) = 2^n - 1.\] 
When $n$ is odd, there are no almost-$\mathrm{SEMI}$-admissible partial orientations. When $n$ is even, there are $2\binom{n}{n/2}$ almost-$\mathrm{SEMI}$-admissible partial orientations: choose half of the edges of $C_n$ and orient them all the same way. Regardless of the parity of $n$, there are $2n$ almost-$\mathrm{SHI}$-admissible partial orientations; there are four cases:
\begin{enumerate}
\item Orient $\{v_1,v_n\}$ as $(v_n,v_1)$ and leave one of the edges $\{v_i,v_{i+1}\}$ blank for~\mbox{$1 \leq i \leq n-1$} while orienting the rest as $(v_i,v_{i+1})$.
\item Orient $\{v_1,v_n\}$ as $(v_1,v_n)$ and leave all the other edges blank.
\item Leave $\{v_1,v_n\}$ blank and orient one of the $\{v_i,v_{i+1}\}$ for $1 \leq i \leq n-1$ as~$(v_{i+1},v_i)$ while leaving leaving the rest blank.
\item Leave $\{v_1,v_n\}$ blank and orient all of the $\{v_i,v_{i+1}\}$ for $1 \leq i \leq n-1$ as~$(v_i,v_{i+1})$. $\square$
\end{enumerate}
\end{example}

\begin{cor}
The maximum number of regions of $\Sigma_G(A)$ over all parameter matrices $A$ is $2^{n-1}T_G(3/2,1)$. This maximum is achieved by $\Sigma_G(\mathrm{GEN})$.
\end{cor}

\emph{Proof}: This follows from Theorem~\ref{thm:genreg} and Proposition~\ref{prop:zeroreg}. $\square$

We can slightly refine the lower bound for the number of regions given by Corollary~\ref{cor:minreg} by considering the degrees of the parking functions of $G_{\sbullet}$.
\begin{definition}
The \emph{degree} of a parking function $c = \sum_{i=1}^{n} c_i v_i$ of $G_{\sbullet}$ with respect to $v_0$ is 
\[\mathrm{deg}(c) := \sum_{i=1}^{n} c_i.\]
\end{definition}

\begin{prop}
Let $g := |E_{\sbullet}| - |V_{\sbullet}| + 1 = |E|$ be the genus of $G_{\sbullet}$. Define the vector~\mbox{$h(G_{\sbullet}) := (h_0, h_1, \ldots, h_g)$}, where 
\[h_i := |\{c\colon c \textrm{ is a parking function of $G_{\sbullet}$ with respect to $v_0$ and } \mathrm{deg}(c) = i\}|.\] 
Define $p(G,A) := (p_0, p_1, \ldots, p_g)$, where 
\[p_i := |\{\mathcal{O}\colon \mathcal{O} \textrm{ is an $A$-admissible partial orientation of $G$ and } |\mathcal{O}| = i\}|.\] 
Then $h(G_{\sbullet}) \leq p(G,A)$.
\end{prop}
\emph{Proof}: This is immediate from the proof of Corollary~\ref{cor:main}: for any region~$R$ of~$\Sigma_G(A)$, we have $\mathrm{deg}(\lambda(R)) = |\mathcal{O}_R|$. $\square$

The lower bound from Corollary~\ref{cor:minreg} is sometimes sharp, as in the case of the complete graph $K_n$.

\begin{prop}
The minimum number of regions of $\Sigma_{K_n}(A)$ over all parameter matrices $A$ is $(n+1)^{(n-1)}$. This minimum is achieved by $\Sigma_{K_n}(\mathrm{SHI})$.
\end{prop}

\emph{Proof}: Shi~\cite{shi} proved that the number of regions of the Shi arrangement is given by~\mbox{$(n+1)^{(n-1)}$}, the number of spanning trees of $K_{n+1}$. $\square$

\begin{remark}
Our Corollary~\ref{cor:main} provides an alternative proof of the bijection of Stanley and Pak~\cite{pak-stanley} between regions of the Shi arrangement and parking functions.
\end{remark}

However, the lower bound from Corollary~\ref{cor:minreg} is in general not sharp:

\begin{example}
Consider the path graph $P_3$, labeled as below:
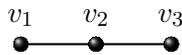
\begin{center}
\begin{tikzpicture}
\SetVertexMath
\GraphInit[vstyle=Art]
\tikzset{VertexStyle/.style = {%
shape = circle,
shading = ball,
ball color = black,
inner sep = 2pt
}}
\SetUpVertex[MinSize=3pt]
\SetVertexLabel
\SetUpEdge[color=black]
\Vertex[LabelOut,Lpos=90,Ldist=.01cm,x=0,y=0]{v_1}
\Vertex[LabelOut,Lpos=90,Ldist=.01cm,x=1,y=0]{v_2}
\Vertex[LabelOut,Lpos=90,Ldist=.01cm,x=2,y=0]{v_3}
\Edge(v_1)(v_2)
\Edge(v_2)(v_3)
\end{tikzpicture}
\end{center}
Then the essentialization of any $\Sigma_{P_3}(A)$ looks like
\begin{center}
\begin{tikzpicture}
\draw [style=thick] (0,0) -- (3,3);
\draw [style=thick] (-2,0) -- (1,3);
\node at (-2.5,-0.5) [text centered,text width=3cm] {$x_3 - x_2 = a_{32}$};
\node at (0.2,-0.5) [text centered,text width=3cm] {$x_2 - x_3 = a_{23}$};
\draw [style=thick] (0,3) -- (3,0);
\draw [style=thick] (-2,3) -- (1,0);
\node at (-2.5,3.5) [text centered,text width=3cm] {$x_1 - x_2 = a_{12}$};
\node at (0.2,3.5) [text centered,text width=3cm] {$x_2 - x_1 = a_{21}$};
\end{tikzpicture}
\end{center}
So $r(\Sigma_{P_3}(A)) = 9$ for any $A$, but there are only $8$ parking functions of ${P_3}_{\sbullet}$ with respect to $v_0$:  $0$, $v_1$, $v_2$, $v_3$, $v_1+v_2$, $v_1 + v_3$, $2v_2$, and $v_2 + v_3$. $\square$
\end{example}

Also, $\Sigma_{G}(\mathrm{SHI})$ does not in general achieve the minimum number of regions of~$\Sigma_G(A)$:

\begin{example}
Consider the cycle graph $C_4$. As was shown in Example~\ref{ex:semicyc},  we have $r(\Sigma_{C_4}(\mathrm{SHI})) = 61$ while $r(\Sigma_{C_4}(\mathrm{SEMI})) = 59$. $\square$
\end{example}

Computing the minimum number of regions of $\Sigma_G(A)$ for arbitrary $G$ and finding a parameter list $A$ that achieves this minimum remain open problems. Proposition~\ref{prop:zeroreg} suggests that maximizing the number of almost-$A$-admissible partial orientations may minimize the number of regions of~$\Sigma_G(A)$, but the exact relationship between almost-admissible orientations and the number of regions remains unclear.

\bibliography{bigraph}{}

\begin{thebibliography}{10}

\bibitem{armstrong}
Drew Armstrong and Brendon Rhoades.
\newblock The {S}hi arrangement and the {I}sh arrangement.
\newblock {\em Trans. Amer. Math. Soc.}, 364(3):1509--1528, 2012.

\bibitem{benson}
Brian Benson, Deeparnab Chakrabarty, and Prasad Tetali.
\newblock {$G$}-parking functions, acyclic orientations and spanning trees.
\newblock {\em Discrete Math.}, 310(8):1340--1353, 2010.

\bibitem{chebikin}
Denis Chebikin and Pavlo Pylyavskyy.
\newblock A family of bijections between {$G$}-parking functions and spanning
  trees.
\newblock {\em J. Comb. Theory Ser. A}, 110(1):31--41, 2005.

\bibitem{duval}
Art Duval, Caroline Klivans, and Jeremy Martin.
\newblock The {$G$}-shi arrangement, and its relation to {$G$}-parking
  functions.
\newblock \url{http://www.math.utep.edu/Faculty/duval/papers/nola.pdf}, January
  2011.

\bibitem{hopkins}
Sam Hopkins and David Perkinson.
\newblock Orientations, semiorders, arrangements, and parking functions.
\newblock {\em Elec. J. of Combin.}, 19(4), 2012.

\bibitem{postnikov}
Alexander Postnikov and Boris Shapiro.
\newblock Trees, parking functions, syzygies, and deformations of monomial
  ideals.
\newblock {\em Trans. Amer. Math. Soc.}, 356(8):3109--3142 (electronic), 2004.

\bibitem{shi}
Jian~Yi Shi.
\newblock {\em The {K}azhdan-{L}usztig cells in certain affine {W}eyl groups},
  volume 1179 of {\em Lecture Notes in Mathematics}.
\newblock Springer-Verlag, Berlin, 1986.

\bibitem{pak-stanley}
Richard~P. Stanley.
\newblock Hyperplane arrangements, interval orders, and trees.
\newblock {\em Proc. Nat. Acad. Sci.}, 93(6):2620--2625, 1996.

\bibitem{stanley}
Richard~P. Stanley.
\newblock An introduction to hyperplane arrangements.
\newblock In {\em Geometric combinatorics}, volume~13 of {\em IAS/Park City
  Math. Ser.}, pages 389--496. Amer. Math. Soc., Providence, RI, 2007.

\bibitem{welsh}
Dominic J.~A. Welsh and Criel Merino.
\newblock The {P}otts model and the {T}utte polynomial.
\newblock {\em J. of Math. Phys.}, 41(3):1127--1152, 2000.

\bibitem{zaslavsky}
Thomas Zaslavsky.
\newblock Facing up to arrangements: face-count formulas for partitions of
  space by hyperplanes.
\newblock {\em Mem. Amer. Math. Soc.}, 1(154):102, 1975.

\end{thebibliography}
\bibliographystyle{plain}

\end{document}